\documentclass[aop,MSNbibl,citesort,dvips]{arximspdf}

\doi{10.1214/11-AOP685}
\volume{40}
\issue{6}
\pubyear{2012}
\firstpage{2439}
\lastpage{2459}

\makeatletter

\newcommand{\rr}{\mathbb{R}}
\newcommand{\rrr}{\mathbb{R}}

\newcommand{\nn}{\mathbb{N}}
\newcommand{\nnn}{\mathbb{N}}

\newcommand{\ee}{\mathbb{E}}

\newcommand{\dd}{\mathbb{D}}

\newcommand{\eu}{\mathrm{e}}
\newcommand{\LL}{\mathrm{L}}

\newtheorem{Proposition}{Proposition}
\newtheorem{Lemma}[Proposition]{Lemma}
\newproclaim{Definition}[Proposition]{Definition}
\newtheorem{Theorem}[Proposition]{Theorem}
\newtheorem{Corollary}[Proposition]{Corollary}

\makeatother

\begin{document}
\begin{frontmatter}

\title{Chaos of a Markov operator and the fourth moment condition}
\runtitle{Chaos of a Markov operator}

\begin{aug}
\author[A]{\fnms{M.} \snm{Ledoux}\corref{}\ead[label=e1]{ledoux@math.univ-toulouse.fr}}
\runauthor{M. Ledoux}
\affiliation{Universit\'{e} de Toulouse and Institut Universitaire de France}
\address[A]{Institut de Math\'{e}matiques de Toulouse\\
Universit\'{e} de Toulouse\\
F-31062 Toulouse\\
France\\
\printead{e1}} 
\end{aug}

\received{\smonth{1} \syear{2011}}
\revised{\smonth{5} \syear{2011}}

%
\begin{abstract}
We analyze from the viewpoint of an abstract Markov operator recent
results by Nualart and Peccati, and Nourdin and Peccati, on
the fourth moment as a condition on a Wiener chaos to have a
distribution close to Gaussian. In particular, we are led to introduce
a notion of chaos associated to a Markov operator through its iterated
gradients and present conditions on the (pure) point spectrum for a
sequence of chaos eigenfunctions to converge to a Gaussian
distribution. Convergence to gamma distributions may be examined
similarly.
\end{abstract}

%

\end{frontmatter}

\section{Introduction}\label{sec1}

In a striking contribution~\cite{Nu-P}, Nualart and Peccati discovered
a few
years ago that the fourth moment
of homogeneous polynomial chaos on Wiener space characterizes
convergence toward
the Gaussian distribution. Specifically, and in a simplified (finite
dimensional) setting, let
$F\dvtx\rr^N \to\rr$, $1 \leq k \leq N$, be defined by
%
%
\begin{equation}\label{equ1}
F = F(x) = \sum_{i_1, \ldots, i_k = 1}^N a_{i_1, \ldots, i_k}
x_{i_1} \cdots x_{i_k},\qquad
x = (x_1, \ldots, x_N) \in\rr^N,
\end{equation}
where $a_{i_1, \ldots, i_k}$ are real numbers vanishing on diagonals
and symmetric in the indices. Assume
by homogeneity that $\int_{\rrr^N} F^2 \,d\gamma_N = 1$ where
\[
d\gamma_N (x) = (2\pi)^{-N/2} \eu^{-|x|^2/2} \,dx
\]
is the standard
Gaussian measure on $\rr^N$.
Such a function $F$ will be called homogeneous of degree $k$. Let now
$F_n$ on $\rr^{N_n}$,
$ n \in\nn$, $N_n \to\infty$, be
a sequence of such homogeneous polynomials of fixed degree $k$. The
main theorem of
Nualart and Peccati~\cite{Nu-P} expresses that the sequence of distributions
of the
$F_n$'s converges\vadjust{\goodbreak} toward the standard Gaussian distribution $\gamma_1$
on the real line if and only if
%
%
\begin{equation}\label{equ2}
\int_{\rrr^{N_n}} F_n^4 \,d\gamma_{N_n} \to3
\end{equation}
(3 being the fourth moment of the standard normal).
The result actually holds for homogeneous chaos on the infinite dimensional
Wiener space, and the equivalence
is further described in terms of convergence of contractions.
The proof of~\cite{Nu-P} relies on multiplication formulas for
homogeneous chaos and the use of stochastic calculus.

Since~\cite{Nu-P} was published, numerous improvements and developments
on this theme have been considered; cf., for example,
\cite{P-Tu,N-OL,No-P1,No-P2,N-P-Rei1,N-P-Rev}$,\ldots.$
An introduction to some of these developments
(with emphasis on multiplication formulas) is the recent monograph
\cite{P-Ta} by Peccati and Taqqu. In particular, the work
by Nualart and Ortiz-Latorre~\cite{N-OL} introduces a technological
breakthrough with a new proof only based on Malliavin calculus and the
use of integration
by parts on Wiener space. In this work, the convergence of ${(F_n)}_{n
\in\nnn}$ to a Gaussian
distribution [and thus also~(\ref{equ2})] is also shown to be
equivalent to the
fact that
%
%
\begin{equation}\label{equ3}
\operatorname{Var}_{\gamma_{N_n}} ( |\nabla F_n |^2) \to0,
\end{equation}
where $\operatorname{Var}_{\gamma_{N_n}}$ is the variance with respect to
$\gamma_{N_n} $.
Based upon this observation, recent work by Nourdin and Peccati
\cite{No-P1,No-P2}
develops the tool of the
so-called Stein method (cf., e.g.,~\cite{S,C-G-S,C-S,R}) in order
to quantify the convergence
toward the Gaussian distribution.
Relying also on multiplication formulas and the use of integration by
parts on Wiener space,
one key step in the investigation~\cite{No-P1} is expressed by the following
inequality:
for a given homogeneous function $F$ of degree $k$ on $\rr^N$
normalized in $\LL^2(\gamma_N)$,
%
%
\begin{equation}\label{equ4}
\operatorname{Var}_{\gamma_N} ( |\nabla F|^2 )
\leq C_k \biggl( \int_{\rrr^N} F^4 \,d\gamma_N - 3 \biggr),
\end{equation}
where $C_k >0$ only depends on $k$.
In particular, the proximity of $ \int_{\rrr^N} F^4 \,d\gamma_N $ to 3
controls the variance
of $ |\nabla F|^2 $. Now, Stein's method for homogeneous
chaos on Wiener space as developed in~\cite{No-P1} expresses that
%
%
\begin{equation}\label{equ5}
d( \nu, \gamma_1) \leq
C \operatorname{Var}_{\gamma_N} ( |\nabla F|^2 ) ^{1/2},
\end{equation}
where $d(\nu, \gamma_1)$ stands for some appropriate distance between
the law $\nu$ of $F$ and $\gamma_1$,
so that $ |\nabla F|^2 $ being close to a constant forces the distribution
of $F$ to be close to a Gaussian distribution. The conjunction of~(\ref{equ4})
and~(\ref{equ5}) thus describes
how the fourth moment condition controls convergence to a Gaussian.

The primary motivation of this work is to understand what
structure of a functional $F$ allows for the preceding results, in
particular thus the control by
the fourth moment of the distance to the Gaussian distribution.
In the process of this investigation, we will revisit\vadjust{\goodbreak}
the preceding results and conclusions in the setting of a symmetric
Markov operator, including, as a particular example, the
Ornstein--Uhlenbeck operator $\LL= \Delta- x \cdot\nabla$,
corresponding to the Wiener space setting.
In order to achieve this goal, observe that
the homogeneous polynomial $F$ of~(\ref{equ1}) is an eigenfunction with
eigenvalue $k$
of the Ornstein--Uhlenbeck operator, that is, $-\LL F = k F$. We shall
therefore try
to understand what is necessary for an eigenfunction $F$ of a Markov
operator in order
to satisfy an inequality such as~(\ref{equ4}). This investigation leads
us to define a notion
of chaos eigenfunction with respect to such a Markov operator with pure
point spectrum
consisting of a countable sequence of eigenvalues,
the homogeneous polynomial $F$ of~(\ref{equ1}) being one example with
respect to the
Ornstein--Uhlenbeck operator. The main achievement
of this work is then the formulation of an explicit
condition on the sequence of eigenvalues under which a chaos
eigenfunction satisfies an inequality such as~(\ref{equ4}).

The basic data will thus be a Markov operator $\LL$ on some state
space $ (E, \mathcal{F}) $ with invariant and reversible probability
measure $\mu$
and symmetric bilinear carr\'{e} du champ operator
\[
\Gamma(f,g) = \tfrac{1}{2} [ \LL(fg) - f \LL g - g \LL f
],
\]
acting on functions $f,g$ in a suitable domain $\mathcal{A}$.
For simplicity, we often write $ \Gamma(f) = \Gamma(f,f)$
which is always nonnegative.
By invariance and symmetry of $\mu$ with respect to $\LL$,
the definition of the carr\'{e} du champ operator
$\Gamma$ yields the integration by parts formula
\[
\int_E f (- \LL g) \,d\mu= \int_E g(- \LL f) \,d\mu= \int_E
\Gamma(f,g) \,d\mu.
\]
In particular $\int_E \LL f \,d\mu= 0$ since $\LL1 = 0$ by the Markov property.
The operator $\LL$ is said, in addition, to be a diffusion operator
if, for every smooth
function $\varphi\dvtx\rr\to\rr$, and every $ f \in\mathcal{A}$,
\[
\LL\varphi(f) = \varphi'(f) \LL f + \varphi'' (f) \Gamma(f).
\]
Alternatively, $\Gamma$ is a derivation in the sense
that $\Gamma(\varphi(f), g) = \varphi'(f) \Gamma(f,g)$.

We refer to the lecture notes~\cite{B}, Chapter 2, by Bakry for an
introduction to this abstract framework
of Markov and carr\'{e} du champ operators and a discussion of some of the
examples emphasized below. Additional general references
include~\cite{D-M-M} for further probabilistic interpretations and
\cite{B-H,F-O-T} for
constructions in terms of Dirichlet forms; see also~\cite{L2} and the
forthcoming~\cite{B-G-L}.
One prototype example of a Markov diffusion operator is the
Ornstein--Uhlenbeck operator acting on say the algebra $\mathcal A$ of
polynomial functions $f$ on
$E = \rr^N$ as $\LL f(x) = \Delta f (x) - x \cdot\nabla f(x)$, with invariant
and reversible probability measure the Gaussian distribution $\mu=
\gamma_N$ and
carr\'{e} du champ $ \Gamma(f) = |\nabla f|^2$. One could consider
its infinite dimensional extension on Wiener space
(cf.~\cite{B-H} and~\cite{N}, Chapter~1), but for simplicity in the\vadjust{\goodbreak}
exposition we
stick here
on the finite dimensional case as a reference example. The preceding
general setting also includes
discrete examples, such as the two-point space and its products.
Namely, on
$ E = \{-1, +1\}^N$, let $ \LL f = {1\over2} \sum_{i=1}^N D_i f $
where $D_i f( x) = f(\tau_i(x)) - f(x)$,
$x = (x_1, \ldots, x_i, \ldots, x_N)$, $\tau_i(x) = (x_1, \ldots, -x_i,
\ldots, x_N)$.
$\LL$ is invariant and symmetric with respect to the uniform measure
$\mu$ on
$ \{-1, +1\}^N$ with carr\'{e} du champ
$ \Gamma(f) = {1\over4} \sum_{i=1}^N (D_i f)^2 $,
but is not a diffusion operator.

These two examples actually entail a crucial
chaos structure in the sense that the generators $\LL$ may be
diagonalized in a
sequence of orthogonal polynomials (Hermite polynomials in the
Gaussian case, Walsh polynomials in the cube example); see, for example,
\cite{B}, Chapter 1,~\cite{N}, Chapter 1,~\cite{J},
Chapter 2,~\cite{P-Ta}, Chapter 5. More
precisely,
setting for $ {\underline k} = (k_1, \ldots, k_N) \in\nn^N$, $x =
(x_1, \ldots, x_N) \in\rr^N$,
$ H _ {{\underline k}} ( x) = h_{k_1} (x_1) \cdots h_{k_N} (x_N) $,
with ${(h_k)}_{k \in\nnn}$ the sequence of orthonormal Hermite
polynomials on the real line,
any function $f\dvtx\rr^N \to\rr$ in $\LL^2(\gamma_N )$ may be written as
\[
f = \sum_{k \in\nnn} \sum_{| {\underline k}| = k}
\langle f , H_{{\underline k}} \rangle H_{{\underline k}},
\]
where $ \langle\cdot, \cdot\rangle$ is the scalar product in
$\LL^2(\gamma_N)$ and
where the second sum runs over all $ {\underline k} \in\nn^N$ with
$ | {\underline k}| = k_1 + \cdots+ k_N = k$.
An element $ H = H_{{\underline k}}$ with $ | {\underline k}|= k$ is an
eigenfunction
of the Ornstein--Uhlenbeck operator with $ - \LL H = k H$ and the spectrum
of the operator $-\LL$ thus consists of the sequence of the
nonnegative integers.
For fixed $k \in\nn$, linear combinations
%
%
\begin{equation}\label{equ6}
F = \sum_{ | {\underline k}| = k} a_ {{\underline k}} H_{{\underline
k}}
\end{equation}
define generic eigenfunctions (chaos) of $-\LL$ with eigenvalue $k$,
the homogeneous
function $F$ of~(\ref{equ1}) being one example.

Similarly, if $f\dvtx\{-1, +1\}^N \to\rr$,
\[
f = \sum_{k = 0}^N \sum_{|A|= k} \langle f , W_A \rangle W_A,
\]
where the second sum runs over all subsets $ A $ of $\{1, \ldots, N\}$
with $k$ elements and
\[
W_A (x) = \prod_{i \in A} x_i,\qquad x = (x_1, \ldots, x_N) \in\{-1,
+1\}^N ,
A \subset\{1, \ldots, N\},
\]
are the\vspace*{1pt} so-called Walsh polynomials. For the discrete operator
$ \LL f =\break {1\over2} \sum_{i=1}^N D_i f $,\vspace*{1pt} $ - \LL W_A = k W_A$ if $|A|=k$.
The spectrum of $-\LL$ is thus equal to $\nn$, and linear combinations
%
%
\begin{equation}\label{equ7}
F = \sum_{|A|=k} a_A W_A
\end{equation}
describe the family of eigenfunctions (chaos) of $-\LL$ with eigenvalue $k$.\vadjust{\goodbreak}

A further example is Poisson space. In dimension one, let
$\mu$ be the Poisson law on $\nn$ with parameter
$\theta>0 $. For a function $ f\dvtx\nn\to\rr$ with finite support
say, let
$Df(j) = f(j) -f(j-1)$ for every $j \in\nn$ [$f(-1) = 0$]. The Poisson operator
may then be defined as
$\LL f(j) = \theta Df(j+1) - j Df(j)$, \mbox{$ j \in\nn$}. It is not a diffusion.
The associated carr\'{e} du champ
operator is given by \mbox{$2 \Gamma(f)(j) = \theta Df(j+1)^2 + j
Df(j)^2$}, $j \in\nn$.
The operator $- \LL$ has a spectrum given by the sequence of the integers
and is diagonalized along the Charlier orthogonal polynomials. Multi-dimensional
Poisson models are similar.

Laplacians $ \LL= \Delta$ on (compact) Riemannian manifolds, and acting
on families of smooth functions, also enter this framework. These Laplacians
are diffusion operators and, in the compact case, have again a spectrum
consisting of a countable sequence of eigenvalues; cf., for example,
\cite{G-H-L}.

This work will analyze properties of eigenfunctions of such Markov
operators $\LL$,
that is, functions $ F\dvtx E \to\rr$ (in the domain of $\LL$) such that
$- \LL F = \lambda F$
for some $\lambda>0$. (We emphasize that $F$ and $\lambda$ are thus rather
eigenfunction and eigenvalue of $-\LL$ which is nonnegative.)
The ultimate goal of this work is to find conditions on such an
eigenfunction $F$ of a diffusion
operator $\LL$ in order that the analog of~(\ref{equ4}) holds, and that the
fourth moment
condition then ensures the proximity with the Gaussian distribution.
We outline here the various steps of the investigation.
The first step will be to show (following~\cite{No-P1} in the
Ornstein--Uhlenbeck
setting) that
Stein's method applied to an eigenfunction $F$ indicates that
it has a Gaussian distribution if (and only if) its
carr\'{e} du champ $\Gamma(F) $ is constant;
see Proposition~\ref{Proposition1} below. More precisely,
in accordance with~(\ref{equ5}), for suitable families of functions
$\varphi\dvtx
\rr\to\rr$,
and whenever $ \int_E F^2 \,d\mu= 1$,
%
%
\begin{equation}\label{equ8}
\biggl| \int_\rrr\varphi(F) \,d\mu- \int_\rrr\varphi\,d\gamma_1
\biggr|
\leq C_\varphi\operatorname{Var}_\mu( \Gamma(F) ) ^{1/2},
\end{equation}
where $\operatorname{Var}_\mu$ is the variance with respect to $\mu$.

On the basis of this result, the fourth moment condition
appears quite naturally by the integration by parts
formula since (assuming the necessary domain and integrability conditions)
\[
\lambda\int_E F^4 \,d\mu=
\int_E F^3 (- \LL F) \,d\mu= 3 \int_E F^2 \Gamma(F) \,d\mu.
\]
Moreover,
$ \int_E \Gamma(F) \,d\mu= \int_E F (-\LL F) \,d\mu= \lambda\int
_E F^2 \,d\mu$,
so that, still assuming by homogeneity that $ \int_E F^2 \,d\mu= 1$,
%
%
\begin{equation}\label{equ9}
\lambda\biggl( {1\over3} \int_E F^4 \,d\mu- 1 \biggr)
= \int_E F^2 \bigl(\Gamma(F) - \lambda\bigr) \,d\mu.
\end{equation}
This identity is the first indication that the proximity of $\int_E F^4
\,d\mu$ with 3 actually
amounts to the proximity of $\Gamma(F) $ with its constant mean value
$\lambda$.

The next step in the investigation, the
main result of this note, describes a chaos structure of an
eigenfunction $F$\vadjust{\goodbreak}
of a Markov operator $\LL$ (not necessarily diffusive) with spectrum
consisting in a sequence
$ S = \{ 0 = \lambda_0 < \lambda_1 < \lambda_2 < \cdots\}$ of eigenvalues
in order that whenever $F$ is such a chaos with
eigenvalue $\lambda_k$ normalized in $\LL^2(\mu)$,
%
%
\begin{equation}\label{equ10}
\operatorname{Var}_\mu( \Gamma(F) )
\leq C_k \int_E F^2 \bigl(\Gamma(F) - \lambda_k \bigr) \,d\mu
\end{equation}
for some finite constant $C_k$ only depending on $S$.
The relations~(\ref{equ8}),~(\ref{equ9}) and~(\ref{equ10}) together
therefore describe how the fourth moment condition\break $\int_E F^4 \,d\mu
\sim3$
ensures that $\Gamma(F)$ is close to constant and thus that
the distribution of $F$ is close to Gaussian.
This family of inequalities may then be used to describe
convergence to a Gaussian distribution of a sequence of such chaos
eigenfunctions.
The abstract chaos structure underlying these results is defined by
means of the iterated gradients of the
Markov operator $\LL$ and is shown to easily cover
the examples of Wiener, Walsh or Poisson chaos. For example, the chaos structure
of the homogeneous polynomial $F$ of~(\ref{equ6}) actually amounts to
the fact that
$\nabla^{k+1} F = 0$. The proof of~(\ref{equ10}) will proceed
by a standard and direct algebraic $\Gamma$-calculus on eigenfunctions
involving the iterated
gradients of the operator $\LL$ and avoiding any type of multiplication
formulas for chaos.

Turning to the content of this note,
Section~\ref{sec2} briefly presents Stein's method applied to an eigenfunction
of a Markov
diffusion operator. The next section discusses the iterated gradients
and the
associated $\Gamma$-calculus on eigenfunctions, of fundamental use
in the investigation. Section~\ref{sec4} introduces the notion of chaos
of a Markov
operator with pure point spectrum
and presents the aforementioned main result~(\ref{equ10}), proved in
Section~\ref{sec6}.
The last section briefly describes analogous conclusions for
convergence to
gamma distributions covering recent results of~\cite{No-P2}.

It should be carefully emphasized that the present exposition develops
more the algebraic and
spectral descriptions of the problem under investigation
[and concentrates on a proof of~(\ref{equ10})]
rather than the analytic
issues on domains and classes of functions involved in the analysis. In
particular, we work
with families of functions in the domain of the Markov operator and its
carr\'{e} du champ and with eigenfunctions assumed to satisfy all the
necessary domain and integrability conditions required to develop
integration by parts and the associated $\Gamma$-calculus. These properties
are classically and easily satisfied for the main examples in mind, the
Gaussian case, the discrete
cube or the setting of the Laplace operator on a compact Riemannian manifold.
Note, however, that the extension from the finite dimensional Gaussian
setting to
the infinite dimensional one requires basic analysis on Wiener space
as presented, for example, in the first chapter of~\cite{N} (see
also~\cite{P-Ta})
in order to fully justify the domain issues and the various
conclusions. These
aspects, carefully developed in the aforementioned references, are not
discussed here.
Further conditions ensuring the validity of the results presented here
might be
developed in broader contexts.\vadjust{\goodbreak}

\section{Stein's method for eigenfunctions}\label{sec2}

We start our investigation with a brief exposition of Stein's lemma applied
to eigenfunctions of a diffusion operator.
We refer to~\cite{S,C-G-S,C-S,R} and the references therein
for general introductions on Stein's method. The results below are mere
adaptations of the investigation~\cite{No-P1} by
Nourdin and Peccati in Wiener space
to which we refer for further details.
Throughout this section, $\LL$ is thus a diffusion operator with
invariant and
reversible measure $\mu$ and carr\'{e} du champ $\Gamma$
as described in the \hyperref[sec1]{Introduction}.
All the necessary domain and integrability conditions on the
eigenfunctions under
investigation are implicitly
assumed, and are satisfied for the main Ornstein--Uhlenbeck
example;
cf.~\cite{No-P1}.

We first illustrate, at a qualitative level, Stein's method in this
abstract context.
Given a measurable map $F\dvtx E \to\rr$, say that $\LL$ commutes to $F$
if there
exists a Markov operator $\mathcal{L}$ on the real line such that for every
$\varphi\dvtx\rr\to\rr$ (in the domain of $\mathcal L$ and such that
$\varphi\circ F$ is in the domain of $\LL$)
\[
\LL(\varphi\circ F) = (\mathcal{L} \varphi) (F).
\]
In this case, the image measure $\mu_F$ of $\mu$ by $F$ is the
invariant measure
of $\mathcal{L}$.

One model factorization operator $\mathcal L$ on $\rr$ is the
Ornstein--Uhlenbeck
operator
$\mathcal{L} \psi= \psi'' - x \psi'$ with invariant measure the standard
Gaussian distribution $d\gamma_1 (x) = \eu^{-x^2/2} \,{dx \over\sqrt
{2\pi}}$.
Let then $ F$ be an eigenfunction of $-\LL$ with eigenvalue \mbox{$\lambda>0$}.
The observation here, at the root of Stein's argument, is that
whenever $ \Gamma= \Gamma(F)$ is ($\mu$-almost everywhere)
constant, then $\LL$ commutes to $F$ through the Ornstein--Uhlenbeck
operator
$\mathcal L$, and thus the distribution $\mu_F$ of $F$ is Gaussian.
Namely, note first
that by integration by parts,\vspace*{1pt}
$ \int_E \Gamma \,d\mu= \int_E F (-\LL F) \,d\mu= \lambda\int_E F^2
\,d\mu$
so that if $\Gamma$ is constant and $F$ is normalized in $\LL^2(\mu
)$, then $\Gamma= \lambda$.
Then, for $\varphi\dvtx\rr\to\rr$ smooth enough, by the chain rule
formula for the diffusion operator $\LL$,
\[
\LL(\varphi\circ F) = \varphi'(F) \LL F + \varphi''(F) \Gamma
= - \lambda F \varphi'(F) + \varphi''(F) \Gamma.
\]
Hence, if $\Gamma= \lambda$,
\[
\LL(\varphi\circ F) = \lambda(\mathcal{L} \varphi) (F)
\]
so that $\LL$ commutes to $F$, and thus $\mu_F$ is the invariant
measure of
the Ornstein--Uhlenbeck operator $\mathcal L$ characterized
as the Gaussian distribution~$\gamma_1 $.

For an eigenfunction $F$, $\Gamma= \Gamma(F)$ constant thus forces
the distribution
of $F$ to be Gaussian. Now, as such, this observation is not of much use
and to describe convergence to normal as for sequences of homogeneous
polynomials
in the \hyperref[sec1]{Introduction}, it should be suitably quantified
in the form of
inequality~(\ref{equ8})
in order to express that the proximity of $\Gamma$ with a constant value
forces the distribution of $F$ to be close to Gaussian. This is the
content of the
classical Stein lemma as described in the next statement.
\begin{Proposition}\label{Proposition1}
Let $ F$ be an eigenfunction of $-\LL$ with
eigenvalue $\lambda>0$
and set $\Gamma= \Gamma(F)$. Denote by $\mu_F$ the distribution of $F$.
Given \mbox{$\varphi\dvtx\rr\to\rr$} integrable with respect to $\mu_F$ and
$\gamma_1$, let\vadjust{\goodbreak}
$\psi$ be a smooth solution of the associated
Stein equation $ \varphi- \int_\rrr\varphi\,d\gamma_1= \psi' - x
\psi$. Then,
%
%
\begin{equation}\label{equ11}
\biggl| \int_\rrr\varphi \,d\mu_F - \int_\rrr\varphi\,d\gamma_1
\biggr|
\leq{ C_\varphi\over\lambda} \biggl(\int_E (\Gamma- \lambda)^2
\,d\mu\biggr)^{1/2},
\end{equation}
where $C_\varphi= {\| \psi'\| }^2_\infty$.
In particular, if $\int_E F^2 \,d\mu= 1$,
\[
\biggl| \int_\rrr\varphi \,d\mu_F - \int_\rrr\varphi\,d\gamma_1
\biggr|
\leq{ C_\varphi\over\lambda} \operatorname{Var}_\mu(\Gamma)^{1/2} .
\]
\end{Proposition}
\begin{pf}
Since $\mu_F$ is the distribution of $F$ under $\mu$, and
by the Stein equation,
\[
\int_\rrr\varphi \,d\mu_F - \int_\rrr\varphi\,d\gamma_1
= \int_E \varphi(F) \,d\mu- \int_\rrr\varphi\,d\gamma_1
= \int_E [ \psi' (F) - F \psi(F) ] \,d\mu.
\]
Now $-\LL F = \lambda F$ so that
\[
\psi' (F) - F \psi(F) = \psi'(F) + \lambda^{-1} \LL F \psi(F)
\]
and hence, after integration by parts with respect to the operator $\LL$
and the use of the diffusion property,
\[
\int_\rrr\varphi \,d\mu_F - \int_\rrr\varphi\,d\gamma_1
= \int_E \psi'(F) [1 - \lambda^{-1} \Gamma] \,d\mu.
\]
Together with the Cauchy--Schwarz inequality,
\[
\biggl| \int_\rrr\varphi \,d\mu_F - \int_\rrr\varphi\,d\gamma_1
\biggr|
\leq\biggl(\int_E \psi'(F)^2 \,d\mu\biggr)^{1/2}
\biggl( \int_E [1 - \lambda^{-1} \Gamma]^2 \,d\mu\biggr)^{1/2},
\]
which amounts to~(\ref{equ11}). If $\int_E F^2 \,d\mu=1$, then
$\int_E \Gamma \,d\mu= \int_E F (-\LL F) \,d\mu= \lambda$
and thus $\int_E (\Gamma- \lambda)^2 \,d\mu= \operatorname{Var}_\mu(\Gamma)$.
The proof of Proposition~\ref{Proposition1} is complete.
\end{pf}

Proposition~\ref{Proposition1} is thus investigated in~\cite{No-P1} for
Wiener chaos.
As is discussed there (Lemma 1.2 and Theorem 3.1), the
constant $C_\varphi$ in~(\ref{equ11}) of Proposition~\ref{Proposition1}
can be uniformly bounded
inside specific
classes of functions.
For instance, $C_\varphi\leq2$ when $\varphi$ is the characteristic
function of a Borel set (corresponding to the total variation distance) and
$C_\varphi\leq1$ when $\varphi$ is the characteristic function of a half-line
(corresponding to the Kolmogorov distance).

For the further purposes, observe, as is classical
(cf.~\cite{S,R}), that Stein's strategy may be developed similarly
for the Laguerre operator on the positive half-line
$\mathcal{L}_p \psi= x \psi'' + (p-x) \psi'$, $ p >0$, with
invariant measure the gamma
distribution $ dg _p(x) = \Gamma(p)^{-1} x^{p-1} \eu^{-x} \,dx $.
Let $ F$ be an eigenfunction of $-\LL$ with eigenvalue\vadjust{\goodbreak} $\lambda>0$
and $\Gamma= \Gamma(F)$.
As above, for every $\varphi\dvtx\rr\to\rr$ smooth enough, setting $ G
= F+p$,
%
\begin{eqnarray*}
\LL(\varphi\circ G)
&=& \varphi'(G) \LL F + \varphi''(G) \Gamma\\
&=& - \lambda F \varphi'(G) + \varphi''(G) \Gamma\\
&=& \lambda\biggl( (p-G) \varphi'(G) + {1\over\lambda} \Gamma
\varphi''(G) \biggr).
\end{eqnarray*}
In this case, if $\Gamma= \lambda G $,
\[
\LL(\varphi\circ G) = \lambda(\mathcal{L}_p \varphi) (G)
\]
so that $\mu_G$ is the invariant measure of $\mathcal{L}_p$ characterized
as the gamma distribution~$g _p$.

For this example of the Laguerre operator, the criterion for an
eigenfunction $F$ to have
a gamma distribution is thus that $\Gamma= \lambda(F+p) $.
On the basis of this qualitative description of Stein's method for the
Laguerre operator, the
next statement illustrates the analog of Proposition~\ref{Proposition1}
for this model.
\begin{Proposition}\label{Proposition2}
Let $ F$ be an eigenfunction of $-\LL$ with
eigenvalue $\lambda>0$, and set $\Gamma= \Gamma(F)$.
Let $p>0$ and denote by $\mu_{F+p}$ the distribution of $F+p$.
Given $\varphi\dvtx\rr\to\rr$ integrable with respect
to $\mu_{F+p}$ and $g_p$, let
$\psi$ be a smooth solution of the associated Stein equation
$ \varphi- \int_\rrr\varphi \, dg_p = x\psi' +(p- x) \psi$. Then,
%
%
\begin{equation}\label{equ12}
\biggl| \int_\rrr\varphi \,d\mu_{F+p} - \int_\rrr\varphi \,dg_p
\biggr|
\leq{C_\varphi\over\lambda}
\biggl( \int_E \bigl( \Gamma- \lambda(F+p) \bigr)^2 \,d\mu\biggr)
^{1/2},
\end{equation}
where $C_\varphi= {\| \psi'\| }^2_\infty$. In particular, if $\int_E
F^2 \,d\mu=p$,
\[
\biggl| \int_\rrr\varphi \,d\mu_{F+p} - \int_\rrr\varphi \,dg_p
\biggr|
\leq{C_\varphi\over\lambda}
\operatorname{Var}_\mu( \Gamma- \lambda F) ^{1/2}.
\]
\end{Proposition}
\begin{pf}
Set again $ G = F+p$. Start as in the proof of Proposition
\ref{Proposition1}, namely
\begin{eqnarray*}
\int_\rrr\varphi \,d\mu_G - \int_\rrr\varphi \,dg_p
&=& \int_E \varphi(G) \,d\mu- \int_\rrr\varphi \,dg_p\\
&=& \int_E [ G \psi' (G) + (p-G) \psi(G) ] \,d\mu.
\end{eqnarray*}
Since $-\LL F = \lambda F$, and thus $\LL G = \lambda( p - G)$,
\[
G \psi' (G) + (p-G) \psi(G) = G \psi' (G) + \lambda^{-1} \LL G
\psi(G) .
\]
After integration by parts with respect to the operator $\LL$ and the
use of the diffusion property,
\[
\int_\rrr\varphi \,d\mu_G - \int_\rrr\varphi \,dg_p
= \int_E \psi'(G) [G - \lambda^{-1} \Gamma] \,d\mu.
\]
The conclusion follows similarly from the Cauchy--Schwarz
inequality.\vadjust{\goodbreak}
\end{pf}

Proposition~\ref{Proposition2} is similarly investigated in \cite
{No-P1} in the context of Stein's method on Wiener space. Again the the
constant $C_\varphi$ in~(\ref{equ12}) may be bounded only in terms of
$p$ inside specific classes of functions; cf.~\cite{No-P1}, Lemma 1.3
and Theorem 3.11. Analogs of Stein's lemma in the context of the
preceding statements have been investigated on discrete Poisson or
Bernoulli spaces in~\cite{P-S-T-U,P-Z,N-P-Rei2}. In those examples, the
control of the variance of $\Gamma$ is not enough to ensure proximity
to a~Gaussian distribution and has to be supplemented by various
additional conditions.

\section{Iterated gradients}\label{sec3}

This section presents the family of the iterated
gradients of a Markov operator and the basic (algebraic) $\Gamma$-calculus
on eigenfunctions at the root of the investigation.
Given a symmetric Markov operator $\LL$ as above (not necessarily a
diffusion operator), recall
following~\cite{B,L1}, the iterated gradients $\Gamma_m$, $ m \geq2$,
associated to
$\LL$ defined according to the rule defining $\Gamma= \Gamma_1$ as
\[
\Gamma_m (f,g) = \tfrac{1}{2} [ \LL\Gamma_{m-1} (f,g) - \Gamma
_{m-1} (f, \LL g)
- \Gamma_{m-1} (g, \LL f) ]
\]
for functions $f, g $ in a suitable class $\mathcal A$.
By extension, $\Gamma_0 (f,g) = fg $.
For simplicity, set $\Gamma_m (f) = \Gamma_m (f,f)$. Note that in general
$\Gamma_m (f)$ for $m \geq2$ is not necessarily nonnegative. The
$\Gamma_2$ operator
has been introduced first by Bakry and \'{E}mery~\cite{B-E} to describe
curvature
properties of Markov operators and to provide a simple criterion to
ensure spectral
gap and functional inequalities; cf.~\cite{B}, Chapter 6,~\cite{L2}
and~\cite{B-G-L}.
This criterion
will be used in Proposition~\ref{Proposition4} below. The iterated
gradients $\Gamma_m$
have been exploited
in~\cite{L1} toward variance and entropy expansions.

The following elementary lemma will be of constant use throughout this note
and concentrates on the significant properties of the iterated
gradients of a given eigenfunction.
Recall that we assume
the necessary domain and integrability conditions to justify the
relevant identities.
\begin{Lemma}\label{Lemma3}
Let $F$ be an eigenfunction of $- \LL$ with
eigenvalue $\lambda$.
Set $\Gamma_m = \Gamma_m(F)$, $m \geq1$.
Then, for every $m \geq1$,
%
%
\begin{equation}\label{equ13}
\Gamma_m = \tfrac{1}{2} \LL\Gamma_{m-1} + \lambda\Gamma_{m-1}
= \bigl(\tfrac{1}{2}\LL+ \lambda\operatorname{Id} \bigr) ^{m-1} \Gamma.
\end{equation}
Furthermore, for every $ m, n \geq1$,
%
%
\begin{equation}\label{equ14}
\int_E \Gamma_n \Gamma_m \,d\mu
= \int_E \Gamma_{n-1} \Gamma_{m +1} \,d\mu.
\end{equation}
In particular, by selecting $n=1$, for every $m \geq1$,
%
%
\begin{equation}\label{equ15}
\int_E \Gamma\Gamma_m \,d\mu= \int_E F^2 \Gamma_{m+1} \,d\mu.
\end{equation}
\end{Lemma}
\begin{pf}
Equality~(\ref{equ13}) is an immediate consequence of the
definition of $\Gamma_m$
and the eigenfunction property
\[
\Gamma_m (F) = \tfrac{1}{2} \LL\Gamma_{m-1} (F) -
\Gamma_{m-1} ( F, \LL F) =
\tfrac{1}{2} \LL\Gamma_{m-1} (F) + \lambda\Gamma_{m-1} (F) .
\]
The conclusion follows by iteration.

Recalling the notation $\Gamma_m = \Gamma_m(F)$, multiply the preceding
identity by $\Gamma_n$ and integrate with respect to $\mu$ to get, by
symmetry,
\[
2 \int_E \Gamma_n \Gamma_m \,d\mu
= \int_E \Gamma_{m-1} \LL\Gamma_n \,d\mu
+ 2 \lambda\int_E \Gamma_ n \Gamma_{m-1} \,d\mu.
\]
Changing the role of $n$ and $m-1$, by symmetry again,
\[
2 \int_E \Gamma_{m-1} \Gamma_{n +1} \,d\mu
= \int_E \Gamma_{m-1} \LL\Gamma_n \,d\mu
+ 2 \lambda\int_E \Gamma_{m-1} \Gamma_n \,d\mu
\]
and the identity~(\ref{equ14}) follows. The proof of the lemma is
complete.
\end{pf}

The following statement is a first illustration of the method developed next.
It expresses a kind of rigidity result under the geometric $\Gamma_2$
curvature condition mentioned previously.
\begin{Proposition}\label{Proposition4}
Assume that the operator $\LL$ is of curvature $\rho>0$ in the sense of
Bakry--\'{E}mery~\cite{B-E} (\cite{B}, Chapter 6), that is, $
\Gamma_2(f) \geq\rho\Gamma(f)$ for every \mbox{$f \in\mathcal{A}$}. If $F$
is an eigenfunction of $-\LL$ with eigenvalue $\rho$, then $\Gamma(F)$
is ($\mu$-almost everywhere) constant. In case $\LL$ is a diffusion
operator, the distribution of $F$ is Gaussian.
\end{Proposition}

It might be useful to recall (\cite{B}, Chapter 6,~\cite{L2,B-G-L})
that under
the curvature condition
of the statement, $\lambda\geq\rho$ for every nonzero
eigenvalue $\lambda$ of $-\LL$.
In particular, $\LL$ is ergodic in the sense that
if $ \Gamma(f) = 0$, then $f$ is constant ($\mu$-almost everywhere).
It is also worthwhile mentioning that for the model space consisting of
the Ornstein--Uhlenbeck
diffusion operator $\LL= \Delta- x \cdot\nabla$ with invariant
measure $\gamma_N $, $\rho=1$
and the eigenfunctions with eigenvalue 1 are the linear functions
\[
F(x) = \sum_{i=1}^N a_i x_i ,\qquad x = (x_1, \ldots, x_N) \in\rr^N,
\]
whose distributions are of course Gaussian. Since Gaussian Wiener
chaos of order larger than or equal to 2 do not contain any nonzero
Gaussian variable~\cite{J} and~\cite{Nu-P}, Proposition~\ref{Proposition4} thus expresses
a kind of
rigidity property
in the sense that if $F$
is a nonzero
eigenfunction of the Ornstein--Uhlenbeck
operator $\LL$ with eigenvalue
$\lambda$, then $F$ is Gaussian if and only if $\Gamma(F)$ is
constant, and if and only
if $ \lambda= \rho=1$.

The proof of Proposition~\ref{Proposition4} is rather straigthforward.
Write as before $\Gamma_m (F) = \Gamma_m$, $m \geq1$. By Lemma~\ref{Lemma3}
[formula~(\ref{equ13})],
$ \Gamma_2 = {1\over2} \LL\Gamma+ \rho\Gamma$.\vadjust{\goodbreak} Therefore,
under the curvature condition $\Gamma_2 (f) \geq\rho\Gamma(f) $,
$ \LL\Gamma\geq0$. But then
\[
0 \leq\int_E \Gamma\LL\Gamma \,d\mu= - \int_E \Gamma( \Gamma
) \,d\mu\leq0 ,
\]
so that $ \Gamma= \Gamma(F)$ is ($\mu$-almost
everywhere) constant. The final assertion of the statement then follows
from Stein's lemma
(Proposition~\ref{Proposition1}).\vspace*{-1pt}

\section{Chaos of a Markov operator}\label{sec4}

This section is devoted to the main conclusions of this work.
We are thus given, on a state space $E$, a Markov operator $\LL$ with
symmetric and
invariant probability measure $\mu$ and carr\'{e} du champ~$\Gamma$
(acting on a suitable algebra of functions $\mathcal A$). Assume in
addition that
$\LL$ has a pure point spectrum consisting of a countable sequence of
eigenvalues
\mbox{$ S = \{ 0 = \lambda_0 < \lambda_1 < \lambda_2 < \cdots\} $}
(more precisely, $S$ is the spectrum of $-\LL$) (cf.~\cite{R-S,Y,B-G-L}).
Since $ \lambda_ 1 >0$, $\LL$ is ergodic
[in the sense that if $ \Gamma(f) = 0$, then $f$ is constant].

Given the spectrum $ S = \{ 0 = \lambda_0 < \lambda_1 < \lambda_2 <
\cdots\} $,
define for every $ k \in\nn$ the polynomial of degree $k$ in the real
variable $X$,
\[
Q_k(X) = \prod_{i=0}^{k-1} (X - \lambda_i)
= \sum_{i=1}^k {1 \over i !} Q_k^{(i)} (0) X^i
\]
($ Q_0 \equiv1)$. Define then the bilinear form (acting on $ \mathcal
{A} \times\mathcal{A}$)
\[
Q_k(\Gamma) = \sum_{i=1}^k {1 \over i !} Q_k^{(i)} (0) \Gamma_i .
\]
The following main definition introduces the notion of chaos associated
to~$\LL$ and its
spectrum $S$.\vspace*{-1pt}
\begin{Definition}\label{Definition5}
An eigenfunction $F$ of $- \LL$ with eigenvalue
$\lambda_k$ ($ - \LL F = \lambda_k F$) is said to be a chaos of
degree $k \geq1$ relative to $S$ if
$Q_{k+1}(\Gamma) (F) = 0$ ($\mu$-almost everywhere).
We call $F$ a chaos eigenfunction (with eigenvalue~$\lambda_k$).\vspace*{-1pt}
\end{Definition}

Motivation for the preceding definition is provided by the Ornstein--Uhlen\-beck
operator
with spectrum $S = \nn$. Namely, it is easily shown in this case
(see~\cite{L1}, Section 2) that
$Q_k(\Gamma) (F) = |\nabla^k F|^2$. Any eigenfunction $F$ as in~(\ref{equ6})
is such that $\nabla^k F$ is constant and $ \nabla^{k+1} F = 0$
leading thus to Definition~\ref{Definition5}.
In the infinite dimensional setting of an abstract Wiener space
$(E, H, \mu)$ with separable Hilbert space $H$, referring
to~\cite{N}, Chapter 1, for notation and terminology,
the Ornstein--Uhlenbeck
operator $\LL$ has domain $\dd^{2,2}$ and
$Q_k(\Gamma)(F) = {\| D^k F\| }^2_{H ^{\otimes k}}$ for any
$ F \in\dd^{k,2}$ where $D $ is the derivative operator
(use as in the finite dimensional case
the commutation $[\LL, D] = D$ and the chain rule formula~\cite{N},
Proposition~1.4.5).
Now, if $J_k F$ denotes the projection of $F$ (in $ \dd^{k,2}$) on the
$k$th Wiener chaos, $\LL J_k F = - k J_k F$ and
$D^k (J_k F) = J_0 D^k F = \ee(D^k F)$ so that $J_kF$ thus
defines a $k$-chaos in the sense of\vadjust{\goodbreak} Definition~\ref{Definition5}.
For example, in case $ H = \LL^2(T, \mathcal{B}, \nu)$ where
$\nu$ is a $\sigma$-finite atomless measure on a measurable space
$(T, \mathcal{B})$, the elements $J_k F$ may be represented as multiple
stochastic integrals
\[
I_k(f_k) = \int_T \cdots\int_T f_k(t_1, \ldots, t_k) W( dt_1)
\cdots W (dt_k)
\]
of symmetric functions $f_k $ on $\LL^2(T^k)$ with respect to the white
noise $W$ and
\[
D^k I_k(f_k) = \{ f_k(t_1, \ldots, t_k) ; t_1, \ldots, t_k
\in T \} .
\]

The discrete operator $ \LL f = {1\over2} \sum_{i=1}^N D_i f $ on the
cube $\{-1,+1\}^N$
and the Poisson operator are further instances entering\vspace*{1pt} this definition
with again $ S = \nn$ (see~\cite{L1}, Section 2). On the cube $\{-1,+1\}^N$,
for example,
\[
Q_k(\Gamma) (F) = {1 \over2^{2k}} \sum(D_{i_1} \cdots D_{i_k} F)^2,
\]
where the sum is over distinct $i_1, \ldots, i_k \in\{1, \ldots, N\}$
and thus
any $ F$ of the form~(\ref{equ7}) is a $k$-chaos ($ k < N$).

There are of course examples of eigenfunctions which are not chaos. For
instance,
the Laguerre operator on the positive half-line $\mathcal{L}_p \psi=
x \psi'' + (p-x) \psi'$, $ p >0$,
has spectrum equal to $\nn$ (with eigenvectors the Laguerre orthogonal
polynomials with respect to the gamma distribution $g_p$), but the
eigenfunction $ F = x - p$ with eigenvalue
1 is not a $1$-chaos as $Q_2(\Gamma) (F) = - {1 \over2} F$.

According to the preceding examples, another possible definition of
$k$-chaos would have been
that $Q_k(\Gamma)(F) $ is constant. (If $F$ is normalized in $\LL
^2(\mu)$, then~(\cite{L1}, page 443),
\[
\int_E Q_k(\Gamma) (F) \,d\mu= \int_E F Q_k (-\LL) F \,d\mu= Q_k
(\lambda_k) ,
\]
hence $Q_k(\Gamma) (F) = Q_k (\lambda_k) $.) Now, it is easily
checked [using~(\ref{equ13}) of Lem\-ma~\ref{Lemma3}]
that if $F$ is an eigenfunction of $- \LL$ with eigenvalue
$\lambda_k$, then $ \LL Q_k(\Gamma)(F) = 2 Q_{k+1}(\Gamma) (F) $.
In particular therefore, if $Q_k(\Gamma)(F)$ is constant, then\break
\mbox{$ Q_{k+1}(\Gamma)(F) = 0$}. Conversely, if $ Q_{k+1}(\Gamma)(F) = 0$, by
ergodicity, $Q_k(\Gamma)(F)$ is constant. It will turn out more simple
in the proofs of the main results to use the first definition of chaos [as
$ Q_{k+1}(\Gamma)(F) = 0$].

The following statements are the main results of this work. Recall
the polynomials $Q_k (X)$ and set, for $k \geq1$, $X \in\rr$,
\[
R_{k+1} (X) = {1 \over X^2} \bigl[ Q_{k+1} (X) - Q_{k+1}^{(1)} (0) X
\bigr]
= \sum_{i=2}^{k+1} {1\over i!} Q_{k+1}^{(i)} (0) X^{i-2}
\]
and
\[
T_{k+1} (X) = R_{k+1} (X + \lambda_k) - R_{k+1} (\lambda_k).
\]
Thus, for example, $Q_2(X) = X^2 - \lambda_1X$,  $R_2 \equiv1$ and $T_2
\equiv0$, $Q_3(X) =  X^3 - (\lambda_1 + \lambda_2) X^2 +\lambda_1
\lambda_2 X$, $R_3(X) = X - (\lambda_1 + \lambda_2)$ and $T_3(X) = X$.
Set furthermore
\[
\pi_k = \lambda_1 \cdots\lambda_k ,\qquad k \geq1 \qquad(\pi_0 =
1).\vadjust{\goodbreak}
\]

The following theorem puts forward the fundamental identity at the root
of this work.\vspace*{-1pt}
\begin{Theorem}\label{Theorem6}
In the preceding setting,
let $F$ be a $k$-chaos eigenfunction
with eigenvalue $\lambda_k$, $k \geq1$. Set $ \Gamma= \Gamma(F)$. Then
%
%
\begin{equation}\label{equ16}
\pi_{k-1} \int_ E \Gamma^2 \,d\mu= \pi_k \int_E F^2 \Gamma \,d\mu
+ (-1)^k \int_E \Gamma T_{k+1} \biggl( { {\LL\over2} }
\biggr) \Gamma \,d\mu.\vspace*{-1pt}
\end{equation}
\end{Theorem}
\begin{Corollary}\label{Corollary7}
In the preceding setting,
let $F$ be a $k$-chaos eigenfunction
with eigenvalue $\lambda_k$, $k \geq1$. Set $ \Gamma= \Gamma(F)$. If
%
%
\begin{equation}\label{equ17}
(-1)^k T_{k+1} \biggl(- { { \lambda_n \over2} } \biggr) \leq0
\qquad\mbox{for every } n \in\nn,
\end{equation}
then
%
%
\begin{equation}\label{equ18}
\int_E \Gamma^2 \,d\mu\leq\lambda_k \int_E F^2 \Gamma \,d\mu.
\end{equation}
In particular, if $F$ is normalized in $\LL^2 (\mu)$, then
$ \int_E \Gamma \,d\mu= \int_E F (- \LL F) \,d\mu= \lambda_k$ and thus
%
%
\begin{equation}\label{equ19}
\operatorname{Var}_\mu(\Gamma) \leq\lambda_k
\biggl(\int_E F^2 \Gamma \,d\mu- \lambda_k \biggr) .\vspace*{-1pt}
\end{equation}
\end{Corollary}

Under the additional diffusion hypothesis on $\LL$, according to~(\ref{equ9}),
inequality~(\ref{equ19}) of Corollary~\ref{Corollary7}
may be expressed equivalently as
%
%
\begin{equation}\label{equ20}
\operatorname{Var}_\mu(\Gamma)
\leq\lambda_k^2 \biggl( {1\over3} \int_E F^4 \,d\mu- 1 \biggr) .
\end{equation}
In particular, if $\int_E F^4 \,d\mu= 3$, then $ \Gamma= \Gamma(F)$
is constant and by Stein's
lemma (Proposition~\ref{Proposition1}), the distribution of $F$ is Gaussian.

The next statement describes a fundamental instance for which the
spectral condition
(\ref{equ17}) in Corollary~\ref{Corollary7} is fulfilled.\vspace*{-1pt}
\begin{Theorem}\label{Theorem8}
The spectral condition~(\ref{equ17}) in Corollary~\ref{Corollary7},
\[
(-1)^k T_{k+1} \biggl(- { { \lambda_n \over2} } \biggr) \leq0
\qquad\mbox{for every } n \in\nn
\]
is satisfied when $S = {(\lambda_n)}_{n \in\nnn} = \nn$.\vspace*{-1pt}
\end{Theorem}

As a consequence of this result, the conclusions of Corollary
\ref{Corollary7} apply to the examples of the Ornstein--Uhlenbeck,
Bernoulli and Poisson operators. As such, some of the main conclusions
of~\cite{No-P1} are covered by the preceding general statement, and in
particular the initial result of~\cite{Nu-P}, namely that if
${(F_n)}_{n \in\nnn}$ is a sequence of homogeneous Gaussian chaos,
normalized in $\LL^2(\gamma_{N_n} )$, $ N_n \to\infty$, then
${(F_n)}_{n \in\nnn}$ converges to a Gaussian distribution as soon as
$\int_E F_n^4 \,d\mu\to3$.\vadjust{\goodbreak}

For discrete models as the cube or the Poisson space, the picture is
less satisfactory.
For instance on the cube $ E = \{-1, +1\}^{N_n}$, $ N_n \to\infty$, if
$ F _n = \sum_{|A|= k} a_A^n W_A $, $ n \in\nn$,
is a sequence of Walsh chaos of degree $k$ normalized in $\LL^2(\mu)$
for the uniform
measure $\mu$, and if $ \int_E F_n^2 \Gamma(F_n) \,d\mu\to k$, then
as an application of Corollary~\ref{Corollary7}, $\Gamma(F_n) \to k$ in
$\LL^2(\mu)$.
Now $\Gamma(F)$ being constant in this case is not always discriminative
[as shown by the example of $F(x) = x_1 \cdots x_k$] and further
conditions have
to be imposed on the sequence ${(F_n)}_{n \in\nnn}$ to ensure convergence
toward a Gaussian distribution. This analysis has been recently
achieved in~\cite{N-P-Rei2}. Similar additional conditions have been
studied on Poisson spaces in~\cite{P-S-T-U,P-Z}. The input of
Corollary~\ref{Corollary7} on
convergence of chaos in these
discrete examples is that it reduces the convergence
$\Gamma(F_n) \to\lambda_k$ in $\LL^2(\mu)$ by the weaker condition
$ \int_E F_n^2 \Gamma(F_n) \,d\mu\to\lambda_k$.

\section{Chaos of order 1 and 2}\label{sec5}

Before turning to the general proofs of Theorem~\ref{Theorem6} and
Corollary~\ref{Corollary7}, and to get a better feeling about these
statements, we discuss in this section the particular values $k=1$ and
$k=2$. Recall that we write for simplicity $\Gamma_m = \Gamma_m(F)$, $m
\geq1$, for an eigenfunction $F$.

When $k=1$, that is, $Q_2(\Gamma) = \Gamma_2 - \lambda_1 \Gamma= 0$,
multiplying this identity~by~$F^2$ and integrating with respect to $\mu$, it follows thanks to
Lemma~\ref{Lemma3} [formula~(\ref{equ15})] that
\[
\int_E \Gamma^2 \,d\mu= \lambda_1 \int_E F^2 \Gamma \,d\mu.
\]
Now here $R_2 \equiv1$, and thus $T_2 \equiv0$, so that both the fundamental
identity~(\ref{equ16}) and the spectral condition~(\ref{equ17}) are
automatically satisfied.

When $k = 2$, start from
$ Q_3(\Gamma) = \Gamma_3 - (\lambda_1 + \lambda_2) \Gamma_2
- \lambda_1 \lambda_2 \Gamma= 0 $.
Multiplying by $F^2$ and integrating, it follows similarly thanks to
Lemma~\ref{Lemma3} [formula~(\ref{equ15})] that
\[
\int_E \Gamma\Gamma_2 \,d\mu- (\lambda_1 + \lambda_2) \int_E
\Gamma^2 \,d\mu
+ \lambda_1 \lambda_2 \int_E F^2 \Gamma \,d\mu= 0 .
\]
By~(\ref{equ13}) of Lemma~\ref{Lemma3},
$\Gamma_2 = {1 \over2} \LL\Gamma+ \lambda_2 \Gamma$ so that
\[
{1 \over2} \int_E \Gamma\LL\Gamma \,d\mu
- \lambda_1 \int_E \Gamma^2 \,d\mu
+ \lambda_1 \lambda_2 \int_E F^2 \Gamma \,d\mu= 0.
\]
Here $R_3(X) = X - (\lambda_1 + \lambda_2)$ and $T_3 (X) = X $ so that
the fundamental identity~(\ref{equ16}) holds, and the spectral
condition~(\ref{equ17})
amounts to $\lambda_n \geq0 $ for every $n \in\nn$.

One observation on which we will come back
in the next section is that, in the case $k = 2$, only the
inequality $ Q_3 (\Gamma) \geq0$ is used in order to reach the
conclusions of Corollary~\ref{Corollary7}.
A further observation is that for chaos of order 1 or 2, the spectral
condition~(\ref{equ17}) is fulfilled
for any sequence of eigenvalues $0 = \lambda_0 < \lambda_1 < \lambda
_2< \cdots.$
This is clearly not the case when $k \geq3$.\vadjust{\goodbreak}

\section{\texorpdfstring{Proofs of Theorems \protect\ref{Theorem6} and \protect\ref{Theorem8}}
{Proofs of Theorems 6 and 8}}\label{sec6}

In this section, we establish Theorem~\ref{Theorem6},
Corollary~\ref{Corollary7} and Theorem~\ref{Theorem8}.
Let thus $F$ be a $k$-chaos with eigenvalue $\lambda_k$.
(If necessary, we may assume that $k \geq3$ according
to the preceding section.) Write as usual
$\Gamma_m$ for $\Gamma_m(F)$, $m\geq1$.

As in the preceding section for chaos of order 1 or 2, start as a first
step from the
chaos hypothesis $ Q_{k+1}( \Gamma) = 0$. Multiply this identity
by $F^2$ and integrate with respect to $\mu$. By definition of $Q_{k+1}$
and~(\ref{equ15}) of Lemma~\ref{Lemma3},
%
%
\begin{eqnarray}\label{equ21}
0 &=& \int_E F^2 Q_{k+1}( \Gamma) \,d\mu\nonumber\\
&=& \sum_{i=1}^{k+1} {1\over i!} Q_{k+1}^{(i)} (0) \int_E F^2
\Gamma_i \,d\mu\\
&=& Q_{k+1} ^{(1)} (0) \int_E F^2 \Gamma \,d\mu
+ \sum_{i=2}^{k+1} {1\over i!} Q_{k+1}^{(i)} (0) \int_E \Gamma
\Gamma_{i-1} \,d\mu. \nonumber
\end{eqnarray}
Now, by~(\ref{equ13}) of Lemma~\ref{Lemma3},
%
\begin{eqnarray*}
&&
\sum_{i=2}^{k+1} {1\over i!} Q_{k+1}^{(i)} (0) \int_E \Gamma
\Gamma_{i-1} \,d\mu\\
&&\qquad= \sum_{i=2}^{k+1} {1\over i!} Q_{k+1}^{(i)} (0)
\int_E \Gamma\biggl( {1\over2} \LL+ \lambda_k \operatorname{Id}
\biggr)^{i-2} \Gamma \,d\mu\\
&&\qquad= \sum_{i=2}^{k+1} {1\over i!} Q_{k+1}^{(i)} (0)
\sum_{\ell= 0}^{i-2} \pmatrix{ {i-2} \cr {\ell} } {1 \over2^\ell}
\lambda_k^{i-2-\ell}
\int_E \Gamma\LL^{\ell} \Gamma \,d\mu\\
&&\qquad= \sum_{\ell= 0}^{k-1} \sum_{i=\ell+ 2}^{k+1}
\pmatrix{ {i-2} \cr
{\ell} } {1\over i!}
Q_{k+1}^{(i)} (0) \lambda_k^{i-2-\ell}
{1 \over2^\ell} \int_E \Gamma\LL^{\ell} \Gamma \,d\mu.
\end{eqnarray*}
Recalling the definition of the polynomial $ R_{k+1}$, note that
\[
\sum_{i=\ell+ 2}^{k+1} \pmatrix{{i-2} \cr{\ell} } {1\over i!}
Q_{k+1}^{(i)} (0) \lambda_k^{i-2-\ell}
= {1 \over\ell!} R_{k+1}^{(\ell)} (\lambda_k) .
\]
Hence
\[
\sum_{i=2}^{k+1} {1\over i!} Q_{k+1}^{(i)} (0) \int_E \Gamma
\Gamma_{i-1} \,d\mu
= \sum_{\ell= 0}^{k-1} {1 \over\ell!} R_{k+1}^{(\ell)} (\lambda
_k) {1 \over2^ \ell}
\int_E \Gamma\LL^\ell\Gamma \,d\mu.
\]
Now
\[
\sum_{\ell= 0}^{k-1} {1 \over\ell!} R_{k+1}^{(\ell)} (\lambda
_k) X ^\ell
= R_{k+1} (X + \lambda_k) = T_{k+1} (X) + R_{k+1} ( \lambda_k)
\]
so that
\[
\sum_{i=2}^{k+1} {1\over i!} Q_{k+1}^{(i)} (0) \int_E \Gamma
\Gamma_{i-1} \,d\mu
= \int_E \Gamma T_{k+1} \biggl( { {\LL\over2} } \biggr)
\Gamma \,d\mu
+ R_{k+1} ( \lambda_k) \int_ E \Gamma^2 \,d\mu.
\]
The fundamental identity~(\ref{equ16}) of Theorem~\ref{Theorem6} then
follows from
(\ref{equ21})
together with the fact that
\[
Q_{k+1}^{(1)} (0) = (-1)^k \lambda_1 \cdots\lambda_k = (-1)^k \pi_k
\]
and
\[
R_{k+1} ( \lambda_k)
= (-1)^{k+1} \lambda_1 \cdots\lambda_{k-1} = (-1)^{k+1} \pi_{k-1}.
\]
The proof is complete.

Corollary~\ref{Corollary7} is deduced from Theorem~\ref{Theorem6}
through the following classical and elementary property, consequence of
the point spectrum hypothesis.
\begin{Lemma}\label{Lemma9}
If $P$ is a polynomial, $\int_E u P(\LL)u \,d\mu
\geq0$ for every
$u$ [in the $\LL^2(\mu)$-domain of $P(\LL)$]
if (and
only if) $P(-\lambda_n) \geq0$ for every $n \in\nn$.
\end{Lemma}
\begin{pf}
For each $n \in\nn$, denote by $E_n$ the eigenspace associated
to the eigenvalue $\lambda_n$ so that $\LL^2 (\mu) = \bigoplus_{n
\in\nnn} E_n$
since $ S = {(\lambda_n)}_{n \in\nnn}$ is the spectrum of~$\LL$.
Decompose then $u$ in $\LL^2(\mu)$ as $ u = \sum_{n \in\nnn} u_n$ with
$ u_n \in E_n$, $ n \in\nn$, so that
$ P (\LL) u = \sum_{n \in\nnn} P( - \lambda_n) u_n$ and
\[
\int_E u P(\LL)u \,d\mu= \sum_{n \in\nnn} P(-\lambda_n) \int_E
u_n^2 \,d\mu
\]
from which conclusion follows.
\end{pf}

As mentioned for chaos of order 2, when $k$ is even, only the inequality
$Q_{k+1}(\Gamma) \geq0$ is used in order to reach the conclusions of
Corollary~\ref{Corollary7}.

We next turn to the proof of Theorem~\ref{Theorem8}, checking the
spectral condition
(\ref{equ17})
$ (-1)^k T_{k+1} (- { \lambda_n \over2} ) \leq0 $, $ n \in
\nn$,
for $S = {(\lambda_n)}_{n \in\nnn} = \nn$.
Since in this case $T_{k+1}(X) = R_{k+1} (X + k) - (-1)^{k+1} (k-1)!$,
we have to show that
\[
\biggl( {n \over2} - k \biggr)^{-2} \Biggl[
\prod_{i=0}^k \biggl( {n \over2} - i \biggr) - k! \biggl( {n \over2} -
k \biggr) \Biggr]
\geq(k-1) ! .
\]
When $ {n \over2} = k$, the expression on the left-hand side
is equal to $ k! \sum_{i=1}^{k} {1\over i}$
so that the conclusion holds in this case.
When $ {n \over2} \not= k$, we need to show that
\[
\biggl( {n \over2} - k \biggr)^{-1} \Biggl[
\prod_{i=0}^{k-1} \biggl( {n \over2} - i \biggr) - k! \Biggr] \geq
(k-1) ! .
\]
Assume first that $n \geq2k +1$. Then
%
\begin{eqnarray*}
\prod_{i=0}^{k-1} \biggl( {n \over2} - i \biggr)
&=& \biggl( {n \over2} - k +1 \biggr) \prod_{i=0}^{k-2} \biggl( {n \over
2} - i \biggr) \\
&\geq& \biggl( {n \over2} - k +1 \biggr) \prod_{i=2}^k \biggl( i + {1
\over2} \biggr)
\geq\biggl( {n \over2} - k +1 \biggr) k ! .
\end{eqnarray*}
Hence
\[
\biggl( {n \over2} - k \biggr)^{-1} \Biggl[
\prod_{i=0}^{k-1} \biggl( {n \over2} - i \biggr) - k! \Biggr] \geq k!,
\]
which answers this case.
We turn to the case where $ n \leq2k -1$ for which it is necessary to
check that
\[
\prod_{i=0}^{k-1} \biggl( {n \over2} - i \biggr) \leq{n \over2}
(k-1)! .
\]
It is enough to assume that $n$ is odd, $n = 2p-1$, $1 \leq p \leq k$. Then
\[
\prod_{i=0}^{k-1} \biggl( {n \over2} - i \biggr)
= \prod_{i=0}^{p-1} \biggl( {n \over2} - i \biggr) \prod_{i=p}^{k-1}
\biggl( {n \over2} - i \biggr)
\leq\prod_{i=1}^p \biggl( i - {1 \over2} \biggr)
\prod_{i=1}^{k-p} \biggl( i - {1 \over2} \biggr) .
\]
Therefore, the inequality to establish amounts to
\[
\prod_{i=1}^{p-1} \biggl( i - {1 \over2} \biggr)
\prod_{i=1}^{k-p} \biggl( i - {1 \over2} \biggr)
\leq(p-1) ! (k-p)! \leq(k-1) !,
\]
which is trivially satisfied. The claims thus holds in this case too.
Theorem~\ref{Theorem8} is therefore established.

\section{Convergence to gamma distributions}\label{sec7}

In this last section, we briefly address the analogs of Theorem
\ref{Theorem6} and Corollary~\ref{Corollary7} in the context of
convergence to gamma distributions on the basis of the corresponding
Stein characterization of Proposition~\ref{Proposition2}. The main
conclusion is obtained by a simple variation on the fundamental
identity~(\ref{equ16}) of Theorem~\ref{Theorem6}. In particular, the
analysis covers the recent results of~\cite{No-P2} (see also
\cite{No-P1}) in the context of Wiener chaos.

The framework is the one of the preceding sections, with a Markov
operator $\LL$
with spectrum $S = {(\lambda_n)}_{n \in\nnn}$ and invariant and reversible
probability measure $\mu$ and carr\'{e} du champ $\Gamma$.
Recall $\pi_k = \lambda_1 \cdots\lambda_k$, $ k \geq1$, and the polynomials
$R_{k+1} $ and $T_{k+1}$ of Theorem~\ref{Theorem6}.

The following theorem addresses approximation of a $k$-chaos $F$ by a
gamma distribution via the
control of $\operatorname{Var}_\mu( \Gamma- \lambda_k F ) $ as emphasized
in Proposition~\ref{Proposition2}.
As announced, the proof is an easy modification on the fundamental
identity~(\ref{equ16}) of Theorem~\ref{Theorem6}.\vadjust{\goodbreak}
\begin{Theorem}\label{Theorem10}
Let $F$ be a $k$-chaos with eigenvalue $\lambda
_k$, $ k \geq1$,
such that $\int_E F^2 \,d\mu=p >0$. Set $ \Gamma= \Gamma(F)$.
Under the spectral condition~(\ref{equ17})\break
$ (-1)^k  T_{k+1} (- { \lambda_n \over2} ) \leq0 $ for
every $n \in\nn$,
it holds
\[
\operatorname{Var}_\mu( \Gamma- \lambda_k F )
\leq\lambda_k \int_E F^2 \Gamma \,d\mu
+ A_k \int_E F \Gamma \,d\mu- p B_k - p^2 \lambda_k^2,
\]
where
\[
A_k = { 2 (-1)^k\lambda_k \over\pi_{k-1} }
R_{k+1} \biggl( {\lambda_k \over2} \biggr)
\quad\mbox{and}\quad
B_k = { (-1)^k\lambda_k^2 \over\pi_{k-1} }
R_{k+1} \biggl( {\lambda_k \over2} \biggr) .
\]
\end{Theorem}

In the diffusion case,
\[
\lambda_k \int_E F^4 \,d\mu= 3 \int_E F^2 \Gamma \,d\mu
\quad\mbox{and}\quad \lambda_k \int_E F^3 \,d\mu= 2 \int_E F
\Gamma \,d\mu
\]
so that the conclusion of the theorem reads
\[
\operatorname{Var}_\mu( \Gamma- \lambda_k F )
\leq{\lambda^2 _k \over3} \int_E F^4 \,d\mu
+ {A_k \lambda_k \over2} \int_E F ^3 \,d\mu- p B_k - p^2 \lambda
_k^2.
\]

Consider now the example where $S= \nn$ for which we know from Theorem~\ref{Theorem8} that
the spectral condition~(\ref{equ17}) holds. The inequality of Theorem
\ref{Theorem10} takes a
nicer form when
$ k \geq2$ is even. Indeed in this case
$ (-1)^k \lambda_k R_{k+1} ( {\lambda_k \over2} ) = - 2 k! $ so that
\[
{1\over k} \operatorname{Var}_\mu( \Gamma- \lambda_k F )
\leq\int_E F^2 \Gamma \,d\mu- 4 \int_E F \Gamma \,d\mu+
2pk - p^2 k .
\]
In particular in the diffusion case,
%
%
\begin{equation}\label{equ22}
{3\over k^2} \operatorname{Var}_\mu( \Gamma- \lambda_k F )
\leq\int_E F^4 \,d\mu- 6 \int_E F^3 \,d\mu+ 6p - 3 p^2 .
\end{equation}
This inequality~(\ref{equ22}) then ensures, through Stein's lemma
(Proposition~\ref{Proposition2}), that if ${(F_n)}_{n \in\nnn}$ is a sequence
of $k$-chaos such that $\int_E F_n ^2 \,d\mu= p$ for every $n$ and
\[
\int_E F_n^4 \,d\mu- 6 \int_E F_n^3 \,d\mu+ 6p - 3 p^2 \to0 ,
\]
then ${(F_n + p )}_{n \in\nnn}$ converges in distribution to
the gamma distribution with parameter $p$, that is the main result of
\cite{No-P2}.
\begin{pf*}{Proof of Theorem~\ref{Theorem10}}
Let thus $F$ be a $k$-chaos with
$\int_E F^2 \,d\mu=p$, hence $\int_E \Gamma \,d\mu= p \lambda_k$.
Set $ U = \Gamma- \lambda_k F$ (so $\int_E U \,d\mu= p\lambda_k$).
It is immediately checked that
%
\begin{eqnarray*}
\int_E \Gamma^2 \,d\mu
&=& \int_E U^2 \,d\mu+ 2 \lambda_k \int_E F \Gamma \,d\mu- p
\lambda_k^2 \\
&=& \operatorname{Var}_\mu(U) + 2 \lambda_k \int_E F \Gamma \,d\mu- p
(1-p) \lambda_k^2
\end{eqnarray*}
and, for every $\ell\geq1$,
\[
\int_E \Gamma\LL^\ell\Gamma \,d\mu
= \int_E U \LL^\ell U \,d\mu+ 2 (-1)^\ell\lambda_k ^{\ell+1}
\int_E F \Gamma \,d\mu
- p (-1)^\ell\lambda_k ^{\ell+2} .
\]
Therefore, the fundamental identity~(\ref{equ16}) of Theorem~\ref{Theorem6}
takes the form, after a little
algebra,
%
\begin{eqnarray*}
&& (-1)^k \int_E U T_{k+1} \biggl( { {\LL\over2} } \biggr) U
\,d\mu
- \pi_{k-1} \operatorname{Var}_\mu(U)\\
&&\qquad{}
+ \pi_k \int_E F^2 \Gamma \,d\mu + 2 (-1)^k \lambda_k R_{k+1} \biggl( { {\lambda_k \over
2} } \biggr)
\int_E F \Gamma \,d\mu\\
&&\qquad{}- p (-1)^k \lambda_k^2 R_{k+1} \biggl( { {\lambda_k \over
2} } \biggr)
- p^2\lambda_k^2 \pi_{k-1} = 0 .
\end{eqnarray*}
Under the spectral condition~(\ref{equ17})
$ (-1)^k T_{k+1} (- { { \lambda_n \over2} } ) \leq0$
for every $n \in\nn$,
%
\begin{eqnarray*}
\pi_{k-1} \operatorname{Var}_\mu(U)
&\leq& \pi_k \int_E F^2 \Gamma \,d\mu
+ 2 (-1)^k\lambda_k R_{k+1} \biggl( { {\lambda_k \over2} }
\biggr)
\int_E F \Gamma \,d\mu\\
&&{} - p (-1)^k \lambda_k^2 R_{k+1} \biggl( { {\lambda_k
\over2} } \biggr)
- p^2 \lambda_k^2 \pi_{k-1},
\end{eqnarray*}
which amounts to the statement of the theorem. The proof is complete.
\end{pf*}

\section*{Acknowledgment}
The author is grateful to one of the referees for the numerous and
constructive comments that led to an improved exposition of this work.


%

%
\printaddresses

\end{document}